\documentclass[a4paper,12pt,reqno]{amsart}
\usepackage{amssymb,amsthm}
\usepackage{ifthen}
\usepackage{graphicx}
\usepackage{float}
\usepackage{mathrsfs}
\usepackage{calligra}
\usepackage[T1]{fontenc}

\theoremstyle{plain}

\newtheorem{thm}{Theorem}[section]
\newtheorem*{thm*}{Theorem}

\numberwithin{equation}{section}

\newtheorem{lem}{Lemma}[section]
\newtheorem{prop}{Proposition}[section]
\newtheorem{rem}{Remark}[section]
\newtheorem{example}{Examples}[section]
\newtheorem{defn}{Definition}[section]

\theoremstyle{definition}
\newcounter {own}
\def\theown {\thesection  .\arabic{own}}

\newenvironment{pf}[1][]{%
 \vskip 3mm
 \noindent
 \ifthenelse{\equal{#1}{}}%
  {{\slshape Proof. }}%
  {{\slshape #1.} }%
 }%
{\qed\bigskip}

\newcounter{alphabet}
\newcounter{tmp}

\newcommand{\ds}{\displaystyle}

\newcounter{minutes}
\divide\time by 60
\newcounter{hours}
\multiply\time by 60 \addtocounter{minutes}{-\time}

\begin{document}
\bibliographystyle{amsplain}
\title{Left translates of a square integrable function on the Heisenberg group}

\thanks{
File:~\jobname .tex,
          printed: \number\year-\number\month-\number\day,
          \thehours.\ifnum\theminutes<10{0}\fi\theminutes}

\author{R. Radha $^\dagger$}

\address{R. Radha, Department of Mathematics,
Indian Institute of Technology Madras, Chennai--600 036, India.}
\email{radharam@iitm.ac.in}
\author{Saswata Adhikari }
\address{Saswata Adhikari, Department of Mathematics,
Indian Institute of Technology Madras, Chennai--600 036, India.}
\email{saswata.adhikari@gmail.com}
\keywords{Besselian, frames, twisted translation, Heisenberg group, Hilbertian.\\
$^\dagger$ {\tt Corresponding author}
}

\maketitle
\pagestyle{myheadings}
\markboth{R. Radha and Saswata Adhikari}{Left translates of a square integrable function on the Heisenberg group}

\begin{abstract}
The aim of this paper is to study some properties of left translates of a square integrable function on the Heisenberg group. First, a necessary and sufficient condition for the existence of the canonical dual to a function $\varphi\in L^{2}(\mathbb{R}^{2n})$ is obtained in the case of twisted shift-invariant spaces. Further, characterizations of $\ell^{2}$-linear independence and the Hilbertian property of the twisted translates of a function $\varphi\in L^{2}(\mathbb{R}^{2n})$ are obtained. Later these results are shown in the case of the Heisenberg group.
\end{abstract}
\section{Introduction}
A closed subspace $V\subset L^{2}(\mathbb{R})$ is called shift-invariant if $f\in V\Longrightarrow \tau_{k}f\in V$ for any $k\in \mathbb{Z}$. Characterizations of shift-invariant spaces in terms of range functions were studied on $\mathbb{R}^{n}$ by Bownik in \cite{bo}. These types of characterization problems were obtained for locally compact abelian groups in \cite{cab,kgr} and for non-abelian compact groups in \cite{rad}. In \cite{bhm}, the authors introduced bracket map on the polarized Heisenberg group $\mathbb{H}_{pol}^{n}$ using the group Fourier transform and obtained characterizations of orthonormal system, frames and Riesz basis consisting of left translates of $\varphi$ in $L^{2}(\mathbb{H}_{pol}^{n})$ in terms of the bracket map. In \cite{cmo}, Currey et al generalized some results of \cite{bo} to shift-invariant spaces associated with a class of nilpotent Lie groups. The concept of the bracket map has been generalized in \cite{bhp} to include any non-abelian discrete group $\Gamma$ using its unitary representations and $L^{1}$ space over the non-commutative measurable space vNa$(\Gamma)$, which is the compact dual of $\Gamma$ whose underlying space is a group von Neumann algebra. Using this bracket map, characterizations of orthonormal basis, Riesz basis, frames were obtained for shift-invariant spaces in a Hilbert space $\mathcal{H}$ given by the action of a non-abelian countable discrete group $\Gamma$. In \cite{lf}, Luef provided a connection between the construction of projections in non-commutative tori and the construction of tight Gabor frames for $L^{2}(\mathbb{R})$. Recently, the authors obtained characterizations of orthonormal system, Bessel sequence, frame and Riesz basis of twisted shift-invariant spaces in terms of the kernel of the Weyl transform in \cite{ras}. Similar characterizations are obtained in the shift-invariant spaces associated with a countably many mutually orthogonal generators on the Heisenberg group in \cite{ra}.

Hernandez et al have provided a necessary and sufficient condition for the existence of the canonical dual to a function $\varphi\in L^{2}(\mathbb{R})$ in \cite{hre}. Further, characterizations of $\ell^{2}$-linear independence and the Hilbertian property of $\{\tau_{k}\varphi:k\in\mathbb{Z}\}$ were obtained in terms of the Fourier transform. The aim of this paper is to obtain similar type of results on the Heisenberg group in terms of the group Fourier transform. In order to obtain our results on the Heisenberg group, we first prove all these results in twisted shift-invariant spaces on $\mathbb{R}^{2n}$. Since $L^{1}(\mathbb{H}^{n})$ is a non-commutative group under convolution and $L^{1}(\mathbb{C}^{n})$ is a non-commutative group under twisted convolution, in order to obtain analogous results, as in the Euclidean case (as per \cite{hre}), we make use of a condition called ``condition C". This condition roughly means that a non-trivial, non-central translate of $\varphi\in L^{2}(\mathbb{H}^{n})$ yields a periodizing (operator valued) sequence on the Fourier transform side that is orthogonal to that of $\varphi$. As mentioned  earlier, we first prove the results in the case of a twisted shift-invariant space on $\mathbb{R}^{2n}$ and then extend to a  shift-invariant space on $\mathbb{H}^{n}$.

The Heisenberg group $\mathbb{H}^n$ is a nilpotent Lie group whose underlying manifold is $\mathbb{R}^n\times\mathbb{R}^{n}\times\mathbb{R}$, where the group operation is defined by $(x,y,t)(x^{\prime},y^{\prime},t^{\prime})=(x+x^{\prime},y+y^{\prime},t+t^{\prime}+\dfrac{1}{2}(x^{\prime}.y-y^{\prime}.x))$ and the Haar measure is the Lebesgue measure  $\mathrm{d}x\mathrm{d}y\mathrm{d}t$ on $\mathbb{R}^n\times\mathbb{R}^{n}\times\mathbb{R}$. Now it is clear that $\{(2k,l,m):k,l\in\mathbb{Z}^{n},m\in\mathbb{Z}\}$ is a discrete subgroup of $\mathbb{H}^{n}$. By Stone-von Neumann theorem, every infinite dimensional irreducible unitary representation on the Heisenberg group is unitarily equivalent to the representation $\pi_{\lambda}, \lambda\in \mathbb{R}^{\star}$, where $\pi_{\lambda}$ is defined by
\begin{eqnarray*}
\pi_{\lambda}(x,y,t)\varphi(\xi)=e^{2\pi i\lambda t}e^{2\pi i\lambda(x.\xi+\frac{1}{2}x.y)}\varphi(\xi+y),
\end{eqnarray*}
where $\varphi\in L^{2}(\mathbb{R}^{n})$.
In order to study shift-invariant spaces on $\mathbb{H}^n$, we need to make use of the representation theory of $\mathbb{H}^n$. The group Fourier transform on $\mathbb{H}^n$ is defined to be $$\widehat{f}(\lambda)=\int\limits_{\mathbb{H}^n}f(x,y,t)\pi_\lambda(x,y,t)dxdydt$$ for $f\in L^{1}(\mathbb{H}^n)$. More explicitly, $\widehat{f}(\lambda)$ is the bounded operator acting on $L^{2}(\mathbb{R}^{n})$ (i.e., $\widehat{f}(\lambda)\in \mathcal{B}(L^{2}(\mathbb{R}^{n}))$) given by $\widehat{f}(\lambda)\varphi=\ds\int\limits_{\mathbb{H}^n}f(x,y,t)\pi_\lambda(x,y,t)\varphi dxdy dt$, where the integral is a Bochner integral taking values in the Hilbert space $L^{2}(\mathbb{R}^{n})$ . Further,
\begin{eqnarray*}
\|\widehat{f}(\lambda)\|_{\mathcal{B}}\leq \|f\|_{L^{1}(\mathbb{H}^{n})}.
\end{eqnarray*}
Define $f^\lambda(x,y)=\ds\int\limits_{\mathbb{R}}f(x,y,t)e^{2\pi i\lambda t}dt$ to be the inverse Fourier transform of $f$ in the $t$-variable. Thus $\widehat{f}(\lambda)=\ds\int\limits_{\mathbb{R}^{2n}}f^\lambda(x,y)\pi_\lambda(x,y,0)dxdy$. One can write $\widehat{f}(\lambda)=W_\lambda(f^\lambda)$, where $W_{\lambda}(f)$ is given by
$$W_\lambda(f)=\int\limits_{\mathbb{R}^{2n}}f(x,y)\pi_\lambda(x,y,0)dxdy,$$
for $f\in L^{1}(\mathbb{R}^{2n})$.

In many problems on $\mathbb{H}^{n}$, an important technique is to take the partial Fourier transform in the $t$-variable to reduce the study to the case of $\mathbb{R}^{2n}$. In particular, for $f,g\in L^{1}(\mathbb{H}^{n})$, the convolution of $f$ and $g$ on $\mathbb{H}^{n}$ is defined to be
\begin{eqnarray*}
(f*g)(z,t)=\int\limits_{\mathbb{H}^{n}}f((z,t)(w,s)^{-1})g(w,s)dwds.
\end{eqnarray*}
This group convolution on $\mathbb{H}^{n}$ can be reduced to $\mathbb{R}^{2n}$ as a non-standard convolution, known as twisted convolution. For $f,g\in L^{1}(\mathbb{R}^{2n})$, the twisted convolution of  $f$ and $g$ is defined to be 
\begin{eqnarray*}
(f\times g)(z)=\int\limits_{\mathbb{R}^{2n}} f(z-w)g(w)e^{\pi i Im (z.\overline{w})}dw.
\end{eqnarray*}
If we define $f^{\#}(z,t)=e^{-2\pi it}f(z)$, then one can show that $f^{\#}* g^{\#}=(f\times g)^{\#}$. Further, for $f,g\in L^{1}(\mathbb{H}^{n})$, one has $\widehat{f*g}(\lambda)=\widehat{f}(\lambda)\widehat{g}(\lambda),~\lambda\in\mathbb{R}^{\star} $, as in the case of Euclidean Fourier transform. This leads to $W(f\times g)= W(f) W(g)$, where 
$$W(f)=\int\limits_{\mathbb{R}^{2n}}f(x,y)\pi(x,y)dxdy,$$ called the Weyl transform of $f\in L^{1}(\mathbb{R}^{2n})$, by taking $\lambda=1$ in $W_{\lambda}(f)$ and by writing $\pi (x,y)=\pi_{1}(x,y,0)$.

Thus in order to study shift-invariant spaces on $\mathbb{H}^{n}$, we consider the twisted shift-invariant spaces on $\mathbb{R}^{2n}$. Let $\mathscr{L}$ be a discrete subgroup of the Heisenberg group $\mathbb{H}^{n}$ such that $\mathbb{H}^{n}/\mathscr{L}$ is compact. In other words, $\mathscr{L}$ is a lattice in $\mathbb{H}^{n}$. For $\varphi\in L^{2}(\mathbb{H}^{n})$, the principal shift-invariant space, $V(\varphi)$, is defined to be  $\overline{span}\{{L_{l}\varphi:l\in\mathscr{L}}\}$, where
$L_{l}\varphi(X)=\varphi(l^{-1}.X),~X\in\mathbb{H}^{n}$. However, for the sake of computational convenience the standard lattice $\{(2k,l,m):k,l\in\mathbb{Z}^{n},m\in\mathbb{Z}\}$ is taken in place of $\mathscr{L}$. Hence for $\varphi\in L^{2}(\mathbb{H}^{n}),~V(\varphi)$ is taken to be the closed linear span of the collection $\{L_{(2k,l,m)}\varphi:k,l\in\mathbb{Z}^{n}, m\in\mathbb{Z}\}$. In order to study the left translations on the Heisenberg group, we consider the twisted translations on $\mathbb{R}^{2n}$. For $\varphi\in L^{2}(\mathbb{R}^{2n}), (k,l)\in\mathbb{Z}^{2n}$, we define the twisted translation of $\varphi$, denoted by $T_{(k,l)}^{t}\varphi$, as
$$T_{(k,l)}^{t}\varphi (x,y)= e^{\pi i(x.l-y.k)}\varphi (x-k,y-l),~(x,y)\in\mathbb{R}^{2n}.$$
Using this definition, the twisted shift-invariant space of $\varphi$, denoted by $V^{t}(\varphi)$, is  defined to be the closed linear span of $\{T_{(k,l)}^{t}\varphi: (k,l)\in\mathbb{R}^{2n}\}$ in $L^{2}(\mathbb{R}^{2n})$.

We shall now mention a few more properties of Weyl transform and group Fourier transform on $\mathbb{H}^{n}$, which will be used in the sequel. The Weyl transform of  a function $f\in L^1(\mathbb{R}^{2n})$ can be explicitly written as
\begin{eqnarray*}
W(f)\varphi (\xi)=\int\limits_{\mathbb{R}^{2n}}  f(x,y)e^{2\pi i(x.\xi+\frac{1}{2}x.y)}\varphi(\xi+y)dxdy,~~ \varphi\in L^{2}(\mathbb {R}^{n}),~~ z=x+iy,
\end{eqnarray*}
which maps $L^{1}(\mathbb {R}^{2n})$ into the space of bounded operators on  $L^{2}(\mathbb {R}^{n})$, denoted by $\mathcal{B}(L^{2}(\mathbb {R}^{n}))$. The Weyl transform $W(f)$ is an integral operator with kernel $K_{f}(\xi,\eta)$ given by
\begin{eqnarray*}
\int\limits_{\mathbb{R}^{n}}f(x,\eta-\xi)e^{i\pi x\cdot(\xi+\eta)}dx.
\end{eqnarray*}
This map $W$ can be uniquely extended to a bijection from the class of tempered distributions $S^{\prime}(\mathbb {R}^{2n})$ onto the space of continuous linear maps from $S(\mathbb {R}^{n})$ into  $S^{\prime}(\mathbb {R}^{n})$. If $f\in L^{2}(\mathbb {R}^{2n})$, then $W(f)\in\mathcal{B}_{2}(L^{2}(\mathbb {R}^{n}))$, the space of Hilbert-Schmidt operators on $\mathbb{R}^{n}$.
For $f,g\in L^{2}(\mathbb {R}^{2n})$, we have
\begin{eqnarray}\label{3peq21}
\langle W(f), W(g)\rangle_{\mathcal{B}_{2}}=\langle f,g\rangle_{L^{2}(\mathbb {R}^{2n})}=\langle K_{f},K_{g}\rangle_{L^{2}(\mathbb {R}^{2n})}.
\end{eqnarray}
The group Fourier transform is an isometric isomorphism of $L^{2}(\mathbb{H}^{n})$ onto $L^{2}(\mathbb{R}^{\star},\mathcal{B}_{2};\\d\mu)$ , where $\mathcal{B}_{2}$ denotes the space of Hilbert-Schmidt operators on $L^{2}(\mathbb{R}^{n})$ and $d\mu(\lambda)=|\lambda|^{n}d\lambda$. For $f,g\in L^{2}(\mathbb{H}^{n})$, we have
\begin{eqnarray}\label{3peq20}
\langle f, g\rangle&=& \int\limits_{\mathbb{R}}\langle \widehat{f}(\lambda),\widehat{g}(\lambda)\rangle_{\mathcal{B}_{2}} |\lambda|^{n}d\lambda
=\int\limits_{\mathbb{R}} \langle W_{\lambda}(f^{\lambda}), W_{\lambda}(g^{\lambda})\rangle_{\mathcal{B}_{2}} |\lambda|^{n}d\lambda.
\end{eqnarray}
We refer to Thangavelu \cite{th} for further details on $\mathbb{H}^{n}$.

The paper is organized as follows. In section 2, we provide required definitions and statement of some results which are available in the literature. In section 3, we study the canonical dual to a function in twisted shift-invariant spaces on $\mathbb{R}^{2n}$. In section 4, we obtain characterization for the twisted translates in $L^{2}(\mathbb{R}^{2n})$ to be $\ell^{2}$-linearly independent. In section 5, we provide characterization for the twisted translates to be Hilbertian. In this case, we show that there exists $\tilde{\varphi}\in L^{2}(\mathbb{R}^{2n})$ such that $\{T_{(k,l)}^{t}\tilde{\varphi}:(k,l)\in\mathbb{Z}^{2n}\}$ is Besselian. In section 6, we obtain these results on the Heisenberg group.
\section{Preliminaries}
Let $\mathcal{H}$ be a separable Hilbert space.
\begin{defn}
A sequence $\{f_{k}:{k \in\mathbb {Z}}\}$ in $\mathcal{H}$ is called a Bessel sequence for $\mathcal{H}$ if there exists a constant $B > 0$ such that
\begin{eqnarray*}
\sum \limits_{k \in\mathbb {Z}}|\langle f,f_{k} \rangle |^{2}\leq B\|f\|^{2},\hspace{4mm}\forall f\in \mathcal{H}.
\end{eqnarray*}
\end{defn}
\begin{defn}
A sequence $\{f_{k}:{k \in\mathbb {Z}}\}$ in $\mathcal{H}$ is called a frame for $\mathcal{H}$ if there exist two constants $A,B > 0$ such that
\begin{eqnarray*}
A\|f\|^{2}\leq\ds\sum \limits_{k \in\mathbb {Z}}|\langle f,f_{k} \rangle |^{2}\leq B\|f\|^{2},\hspace{4mm}\forall f\in \mathcal{H}.
\end{eqnarray*}
\end{defn}
\begin{defn}
A sequence $\{f_{k}:{k \in\mathbb {Z}}\}$ in $\mathcal{H}$ is said to be $\ell^{2}$-linearly dependent if there exists a non-zero sequence $\{c_{k}\}\in\ell^{2}(\mathbb{Z})$ such that $\sum\limits_{k\in\mathbb{Z}}c_{k}f_{k}=0$.
If the sequence $\{c_{k}\}$ is not $\ell^{2}$-linearly dependent, then it is said to be  $\ell^{2}$-linearly independent.
\end{defn}
For a study of of frames we refer to \cite{chr2} and \cite{heil}.

We shall make use of the following definitions and results which were given in \cite{ras}.
\begin{lem}
Let $\varphi \in L^{2}(\mathbb {R}^{2n})$. Then the kernel of the Weyl transform of $T_{(k,l)}^{t}\varphi$ satisfies the following relation.
\begin{eqnarray}\label{3peq1}
K_{T_{(k,l)}^{t}\varphi}(\xi,\eta)=e^{\pi i (2\xi+l).k}K_{\varphi}(\xi+l,\eta).
\end{eqnarray}
\end{lem}
\begin{defn}
For $\varphi\in L^{2}(\mathbb{R}^{2n})$, the function $w_{\varphi}$ is defined as follows.
\begin{eqnarray*}
w_{\varphi}(\xi)=\sum\limits_{m\in \mathbb{Z}^{n}}\int\limits_{\mathbb{R}^{n}}| K_{\varphi}(\xi+m,\eta)|^{2} d\eta,\hspace{4 mm} \xi\in \mathbb{R}^{n}.
\end{eqnarray*}
\end{defn}
\begin{defn}
A function $\varphi \in L^{2}(\mathbb {R}^{2n})$ is said to satisfy ``condition C'' if
\begin{eqnarray*}
\sum\limits_{m\in\mathbb {Z}^{n}}\int\limits_{\mathbb {R}^{n}}K_{\varphi}(\xi+m,\eta)\overline{K_{\varphi}(\xi+m+l,\eta)}d\eta=0 ~~ a.e.~ \xi\in \mathbb {T}^{n}, ~~{\rm for~ all} ~~l\in\mathbb{Z}^{n}\setminus \{0\}.
\end{eqnarray*}
\end{defn}
\begin{thm}\label{3pth10} \cite{ras}
If $\{T_{(k,l)}^{t}\varphi:(k,l)\in\mathbb {Z}^{2n}\}$ is a Bessel sequence in $L^{2}(\mathbb {R}^{2n})$ with bound B, then $w_{\varphi}(\xi)\leq B$ a.e. $\xi\in\mathbb {T}^{n}$. Conversely, suppose  $w_{\varphi}(\xi)\leq B$ a.e. $\xi\in\mathbb {T}^{n}$. If, in addition $\varphi$ satisfies condition C, then $\{T_{(k,l)}^{t}\varphi:(k,l)\in\mathbb {Z}^{2n}\}$ is a Bessel sequence in $L^{2}(\mathbb {R}^{2n})$ with bound B.
\end{thm}

Let $\varphi\in L^{2}(\mathbb{R}^{2n})$ be such that $\varphi$ satisfies condition C. Suppose $A^{t}(\varphi)={\rm span}\{T_{(k,l)}^{t}\varphi\\:{(k,l)\in \mathbb{Z}^{2n}}\}$ and $V^{t}(\varphi)=\overline{A^{t}(\varphi)}$. Consider $f\in A^{t}(\varphi)~i.e., f= \sum\limits_{(k^{\prime},l^{\prime})\in\mathcal{F}}c_{k',l'}T_{(k',l')}^{t}\varphi$, where $\mathcal{F}$ is a finite set. Define $\rho(\xi)=\{\fontdimen16\textfont2=5pt \rho_{l^{\prime}}(\xi)\}_{l^{\prime}\in\mathbb {Z}^{n}}$ for $\xi\in\mathbb {T}^{n}$, where
$\fontdimen16\textfont2=5pt \rho_{l^{\prime}}(\xi)=\sum\limits_{k^{\prime}}c_{k^{\prime},l^{\prime}}e^{\pi i(2\xi+l^{\prime}).k^{\prime}}$. Define $J_{\varphi}(f)=\rho$. In particular, taking $f=T_{(k,l)}^{t}\varphi$, one has
\begin{eqnarray}
J_{\varphi}(T_{(k,l)}^{t}\varphi)(\xi)=(\dotsc ,0,\dotsc ,0,e^{\pi il.k}e^{2\pi ik.\xi},0,\dotsc ,0,\dotsc)\label{3peq8}
\end{eqnarray}
with $e^{\pi il.k}e^{2\pi ik.\xi}$ in the $l$th position a.e. $\xi\in\mathbb{T}^{n}$.
\begin{prop}\label{3pprop1} \cite{ras}
The map $J_{\varphi}$ initially defined on $A^{t}(\varphi)$ can be extended to an isometric isomorphism between $V^{t}(\varphi)$ and $L^{2}(\mathbb {T}^{n},\ell^{2}(\mathbb {Z}^{n}),w_{\varphi})$.
\end{prop}
Moreover, it was proved that  $f\in V^{t}(\varphi)$ if and only if
\begin{eqnarray}\label{3peq2}
K_{f}(\xi,\eta)=\ds\sum\limits_{l^{\prime}\in\mathbb {Z}^{n}}\fontdimen16\textfont2=5pt \rho_{l^{\prime}}(\xi)K_{\varphi}(\xi+l',\eta),
\end{eqnarray}
where $\rho(\xi)= \{\fontdimen16\textfont2=5pt \rho_{l^{\prime}}(\xi)\}_{l^{\prime}\in\mathbb {Z}^{n}}$ and $\rho\in L^{2}(\mathbb {T}^{n},\ell^{2}(\mathbb {Z}^{n}),w_{\varphi})$.

The following equation \eqref{3peq6} appears in the proof of  Theorem 3.5 in \cite{ras}. 
\begin{lem}
Let $\{c_{k,l}:(k,l)\in\mathbb{Z}^{2n}\}$ be a finite sequence and $\varphi\in L^{2}(\mathbb{R}^{2n})$ be such that $\varphi$ satisfies condition C. Then
\begin{eqnarray}
\bigg\|\sum\limits_{(k,l)\in\mathcal{F}}c_{k,l}T_{(k,l)}^{t}\varphi\bigg\|_{L^{2}(\mathbb{R}^{2n})}^{2}=\sum\limits_{l}
\int\limits_{\mathbb {T}^{n}}\bigg|\sum\limits_{k}c_{k,l}e^{\pi il.k}e^{2\pi ik.\xi}\bigg|^{2} w_{\varphi}(\xi)d\xi,\label{3peq6}
\end{eqnarray}
where $\mathcal{F}$ denotes a finite set.
\end{lem}
Now, we shall give some definitions and results which are also given in \cite{ra}.
\begin{defn}
For $\varphi\in L^{2}(\mathbb{H}^{n})$ and $(k,l)\in\mathbb{Z}^{2n}$, the function $G_{k,l}^{\varphi}$ is defined as follows.
\begin{eqnarray}\label{3peq18}
G_{k,l}^{\varphi}(\lambda)=\sum\limits_{r\in\mathbb{Z}}\langle\widehat{\varphi}
(\lambda+r),\widehat{L_{(2k,l,0)}\varphi}(\lambda+r)\rangle_{\mathcal{B}_{2}}
|\lambda+r|^{n},~\lambda\in (0,1].
\end{eqnarray}
\end{defn}
In fact, originally in \cite{ra}, the function $G_{k,l}^{\varphi}$ was defined in terms of the kernel of $W_{\lambda}$. Later it was shown that $G_{k,l}^{\varphi}$ turns out to be \eqref{3peq18}.
\begin{defn}
A function $\varphi\in L^{2}(\mathbb{H}^{n})$ is said to satisfy ``condition C'' if $G_{k,l}^{\varphi}(\lambda)=0$ a.e. $\lambda\in (0,1],~\forall~(k,l)\in\mathbb{Z}^{2n}\setminus\{(0,0)\}$.
\end{defn}
\begin{rem}
In order to show that a function $\varphi\in L^{2}(\mathbb{H}^{n})$ satisfies condition C, it is enough to show that all the Fourier coefficients of $G_{k,l}^{\varphi}$ vanish when $(k,l)\in\mathbb{Z}^{2n}\setminus \{(0,0)\}$.
But
\begin{eqnarray}
\widehat{G_{k,l}^{\varphi}}(m)=\int\limits_{0}^{1} G_{k,l}^{\varphi}(\lambda) e^{-2\pi im\lambda}d\lambda
=\langle\varphi, L_{(2k,l,m)}\varphi\rangle_{L^{2}(\mathbb{H}^{n})}.\label{3pr6}
\end{eqnarray}
Thus it is enough to show that $\langle\varphi, L_{(2k,l,m)}\varphi\rangle=0,~\forall~m\in\mathbb{Z}$, whenever $(k,l)\neq (0,0)$.
\end{rem}
\begin{example}
We shall first provide some examples of functions in $L^{2}(\mathbb{H})$ which satisfy condition C.
\end{example}
\begin{itemize}
\item [1.] Let $\varphi(x,y,t)=\chi_{[0,2]}(x)\chi_{[0,1]}(y)h(t)$, where $h$ is an arbitrary function in $L^{2}(\mathbb{R})$. Then
\begin{eqnarray*}
L_{(2k,l,m)}\varphi (x,y,t)
&=&\varphi(x-2k,y-l,t-m+\frac{1}{2}(2ky-xl))\\
&=& \chi_{[0,2]}(x-2k)\chi_{[0,1]}(y-l) h(t-m+ky-\frac{l}{2}x)\\
&=& \chi_{[2k,2k+2]}(x)\chi_{[l,l+1]}(y) h(t-m+ky-\frac{l}{2}x).
\end{eqnarray*}
Since for $(k,l)\neq (0,0), [0,2]\times [0,1]~\bigcap~ [2k,2k+2]\times [l,l+1]=\emptyset,$  it follows that $\langle \varphi, L_{(2k,l,m)}\varphi\rangle=0,~\forall~ m\in\mathbb{Z}$. Then from \eqref{3pr6}, we get $\widehat{G_{k,l}^{\varphi}}(m)=0,~\forall~ m\in\mathbb{Z}$ whenever $(k,l)\neq (0,0)$. Thus $\varphi$ satisfies condition C.

In \cite{ra}, we have proved the theorem that $\{L_{(2k,l,m)}\varphi:(k,l,m)\in\mathbb{Z}^{2n+1}\}$ is an orthonormal system in $L^{2}(\mathbb{H}^{n})$ if and only if $G_{0,0}^{\varphi}(\lambda)=1$ a.e. $\lambda\in (0,1]$ and $\varphi$ satisfies condition C. Thus it is meaningful to give an example of a function $\varphi$ which satisfies condition C but $\{L_{(2k,l,m)}\varphi:(k,l,m)\in\mathbb{Z}^{2n+1}\}$ is not an orthonormal system. This is illustrated in Example 2 and Example 3.
\item [2.] We take $h(t)=e^{-t^{2}}, t\in\mathbb{R},$ in Example 1. For this function $h$, the correspond-ing $\varphi$ satisfies $\langle\varphi, L_{(0,0,1)}\varphi\rangle\neq 0$.
Hence $\{L_{(2k,l,m)}\varphi:(k,l,m)\in\mathbb{Z}^{2n+1}\}$ does not form an orthonormal system in $L^{2}(\mathbb{H})$.
\item [3.] Instead of $e^{-t^{2}}$, one can take $e^{-|t|}$ or in general any $h\in L^{2}(\mathbb{R})$ for which $\int\limits_{\mathbb{R}}h(t)h(t-1)dt\neq 0$. Then $\langle\varphi, L_{(0,0,1)}\varphi\rangle\neq 0$.
\item [4.] Let $\varphi(x,y,t)=f(x)g(y)h(t)$, where $supp~f\bigcap supp~\tau_{2k}f=\emptyset,~\forall~k\neq 0$ and $supp~g\bigcap supp~\tau_{l}g=\emptyset,~\forall~l\neq 0,~h\in L^{2}(\mathbb{R})$. Then the same proof discussed in Example 1 can be used to show that $\varphi$ satisfies condition C. Again as an example, we can take
\begin{eqnarray*}
f(x)&=\left\{
\begin{array}{l l}
e^{-\frac{1}{x(2-x)}},& \quad \text{$0<x<2$,}\\
0, & \quad \text{otherwise,}
\end{array}\right.
\end{eqnarray*}
and
\begin{eqnarray*}
g(x)&=\left\{
\begin{array}{l l}
e^{-\frac{1}{x(1-x)}},& \quad \text{$0<x<1$,}\\
0, & \quad \text{otherwise,}
\end{array}\right.
\end{eqnarray*}
and $h(t)=e^{-t^{2}},~t\in\mathbb{R}$.

Now we shall provide examples of functions in $L^{2}(\mathbb{H})$, which do not satisfy condition C.
\item [5.]Let $\varphi(x,y,t)=\chi_{[0,3]}(x)\chi_{[0,1]}(y)h(t)$, where $h$ is an arbitrary function in $L^{2}(\mathbb{R})$. Then it can be shown that 
\begin{eqnarray*}
\langle\varphi, L_{(2,0,0)}\varphi\rangle
=\int\limits_{0}^{1}\int\limits_{\mathbb{R}} h(t)\overline{h(t+y)}dtdy.
\end{eqnarray*}
Choose $h$ in such a way that the above integral is non-zero. For example, we can take $h(t)=sinct,~t\in\mathbb{R}$. Then
\begin{eqnarray*}
\langle\varphi, L_{(2,0,0)}\varphi\rangle
=\int\limits_{0}^{1}\frac{sin\pi y}{\pi y}dy
>0.
\end{eqnarray*}
Thus from \eqref{3pr6}, we have $\widehat{G_{1,0}^{\varphi}}(0)\neq 0$, showing that $\varphi$ does not satisfy condition C. 

\item [6.] More generally, let $\varphi(x,y,t)=f(x)g(y)h(t)$, where the value of both the integrals $\int\limits_{\mathbb{R}} f(x)\overline{f(x-2)}dx$ and $\int\limits_{\mathbb{R}}\int\limits_{\mathbb{R}}|g(y)|^{2}h(t)\overline{h(t+y)}dt dy$ are non-zero. Since
$\langle\varphi, L_{(2,0,0)}\varphi\rangle\neq 0$, $\varphi$ does not satisfy condition C. For example, we can take $f(x)=e^{-|x|},~x\in\mathbb{R},~ h(t)=sinct,~t\in\mathbb{R}$ and 
\begin{eqnarray*}
g(y)&=\left\{
\begin{array}{l l}
\frac{1}{n+1},& \quad \text{$y\in [2n,2n+1],~n=0,1,2,\dotsc $,}\\
0, & \quad \text{otherwise.}
\end{array}\right.
\end{eqnarray*}
 \end{itemize}
\begin{rem}
Let $m=0$ in \eqref{3pr6}. Then we have
\begin{eqnarray*}
\int\limits_{0}^{1}G_{k,l}^{\varphi}(\lambda)d\lambda=\langle\varphi, L_{(2k,l,0)}\varphi\rangle<\infty,
\end{eqnarray*}
which shows that the function $G_{k,l}^{\varphi}(\lambda)$, defined in \eqref{3peq18}, is finite a.e. $\lambda\in (0,1]$.
\end{rem}

The following theorem is a consequence of Theorem 4.1 in \cite{ra}. For its proof, we refer to \cite{ra}.
\begin{thm}\label{3pth11}
If $\{L_{(2k,l,m)}\varphi:(k,l,m)\in\mathbb{Z}^{2n+1}\}$ is a Bessel sequence in $L^{2}(\mathbb{H}^{n})$ with bound B, then $G_{0,0}^{\varphi}(\lambda)\leq B$ a.e. $\lambda\in(0,1]$. Conversely, suppose $G_{0,0}^{\varphi}(\lambda)\leq B$ a.e. $\lambda\in(0,1]$. If, in addition $\varphi$ satisfies condition C, then $\{L_{(2k,l,m)}\varphi:(k,l,m)\in\mathbb{Z}^{2n+1}\}$ is a Bessel sequence in $L^{2}(\mathbb{H}^{n})$ with bound B.
\end{thm}

Let $\varphi\in L^{2}(\mathbb{H}^{n})$ be such that $\varphi$ satisfies condition C. Suppose $A(\varphi)=$ span $\{L_{(2k,l,m)}\varphi:(k,l,m)\in\mathbb{Z}^{2n+1}\}$. Then  $V(\varphi)=\overline{A(\varphi)}$. Let $f\in A(\varphi)$ i.e., $f=\sum\limits_{(k^{\prime},l^{\prime},m^{\prime})\in\mathcal{F}}
c_{k^{\prime},l^{\prime},m^{\prime}}L_{(2k^{\prime},l^{\prime},m^{\prime})}\varphi$, where $\mathcal{F}$ is a finite set. Define $\rho(\lambda)=\fontdimen16\textfont2=5pt \{\rho_{k^{\prime},l^{\prime}}(\lambda)\}_
{(k^{\prime},l^{\prime})}$ for $\lambda\in (0,1]$, where
\begin{eqnarray}\label{3pr1}
\fontdimen16\textfont2=5pt \rho_{k^{\prime},l^{\prime}}(\lambda)=\sum\limits_{m^{\prime}}
c_{k^{\prime},l^{\prime},m^{\prime}}e^{2\pi im^{\prime}\lambda}.
\end{eqnarray}
Define $J_{\varphi}(f)=\rho$.
\begin{prop}
The map $J_{\varphi}$ initially defined on $A(\varphi)$ can be extended to an isometric isomorphism between $V(\varphi)$ and $L^{2}((0,1],\ell^{2}(\mathbb{Z}^{2n}),G_{0,0}^{\varphi})$.
\end{prop}
(The above result has been proved for a more general case in \cite{ra}. However, for the sake of completeness, we provide the proof here.)
\begin{pf}
Let $f\in A(\varphi)$. Then
\begin{eqnarray}
\widehat{f}(\lambda)&=&\sum\limits_{(k^{\prime},l^{\prime},m^{\prime})\in\mathcal{F}}
c_{k^{\prime},l^{\prime},m^{\prime}}\widehat {L_{(2k^{\prime},l^{\prime},m^{\prime})}\varphi}(\lambda)\nonumber\\
&=&\sum\limits_{k^{\prime},l^{\prime},m^{\prime}}c_{k^{\prime},l^{\prime},m^{\prime}} e^{2\pi im^{\prime}\lambda}\widehat {L_{(2k^{\prime},l^{\prime},0)}\varphi}(\lambda)\nonumber\\
&=&\sum\limits_{k^{\prime},l^{\prime}}\fontdimen16\textfont2=5pt \rho_{k^{\prime},l^{\prime}}(\lambda)\widehat {L_{(2k^{\prime},l^{\prime},0)}\varphi}(\lambda),\label{3pr5}
\end{eqnarray}
using \eqref{3pr1} and $\rho(\lambda)=\fontdimen16\textfont2=5pt \{\rho_{k^{\prime},l^{\prime}}(\lambda)\}_
{(k^{\prime},l^{\prime})}$.

Conversely let $\rho(\lambda)=\fontdimen16\textfont2=5pt \{\rho_{k^{\prime},l^{\prime}}(\lambda)\}_
{(k^{\prime},l^{\prime})\in\mathcal{F}}$, where $\rho_{k^{\prime},l^{\prime}}(\lambda)$ is given by \eqref{3pr1}. Define $f=\sum\limits_{(k^{\prime},l^{\prime},m^{\prime})\in\mathcal{F}}
c_{k^{\prime},l^{\prime},m^{\prime}}L_{(2k^{\prime},l^{\prime},m^{\prime})}\varphi$. Then $f\in A(\varphi)$. Thus we see that there is a one to one correspondence between $A(\varphi)$ and the collection of functions of the form $\rho$. Further, we have
\begin{eqnarray}
\|f\|^{2}&=&\bigg\|\sum\limits_{(k^{\prime},l^{\prime},m^{\prime})\in\mathcal{F}}
c_{k^{\prime},l^{\prime},m^{\prime}}L_{(2k^{\prime},l^{\prime},m^{\prime})}\varphi\bigg\|^{2}\nonumber\\
&=&\sum\limits_{k^{\prime},l^{\prime}}\int\limits_{0}^{1}\bigg|\sum\limits_{m^{\prime}}c_{k^{\prime},l^{\prime},m^{\prime}}e^{2\pi im^{\prime}\lambda}\bigg|^{2}G_{0,0}^{\varphi}(\lambda)d\lambda.\label{3peq14}
\end{eqnarray}
In fact,
\begin{eqnarray}
&&\bigg\|\sum\limits_{(k^{\prime},l^{\prime},m^{\prime})\in\mathcal{F}} c_{k^{\prime},l^{\prime},m^{\prime}}L_{(2k^{\prime},l^{\prime},m^{\prime})}\varphi\bigg\|_{L^{2}(\mathbb{H}^{n})}^{2}\nonumber\\
&=&\int\limits_{\mathbb{R}}\bigg\|\sum\limits_{k^{\prime},l^{\prime},m^{\prime}}c_{k^{\prime},l^{\prime},m^{\prime}}
\widehat{L_{(2k^{\prime},l^{\prime},m^{\prime})}\varphi(\lambda)}\bigg\|_{\mathcal{B}_{2}}^{2}
|\lambda|^{n}d\lambda\nonumber\\
&=&\int\limits_{0}^{1}\sum\limits_{r\in\mathbb{Z}}\bigg\|\sum\limits_{k^{\prime},l^{\prime},m^{\prime}}
c_{k^{\prime},l^{\prime},m^{\prime}}\widehat{L_{(2k^{\prime},l^{\prime},m^{\prime})}\varphi}(\lambda+r)\bigg\|
_{\mathcal{B}_{2}}^{2}|\lambda+r|^{n}d\lambda\nonumber\\
&=&\int\limits_{0}^{1}\sum\limits_{r\in\mathbb{Z}}\bigg\|\sum\limits_{k^{\prime},l^{\prime}}
\bigg(\sum\limits_{m^{\prime}}c_{k^{\prime},l^{\prime},m^{\prime}} e^{2\pi im^{\prime}\lambda}\bigg)\widehat{L_{(2k^{\prime},l^{\prime},0)}\varphi}(\lambda+r)\bigg\|
_{\mathcal{B}_{2}}^{2}|\lambda+r|^{n}d\lambda\nonumber
\end{eqnarray}
\begin{eqnarray}
&=&\int\limits_{0}^{1}\sum\limits_{r\in\mathbb{Z}}\sum\limits_{k^{\prime},l^{\prime}}\bigg\|
\bigg(\sum\limits_{m^{\prime}}c_{k^{\prime},l^{\prime},m^{\prime}} e^{2\pi im^{\prime}\lambda}\bigg)\widehat{L_{(2k^{\prime},l^{\prime},0)}\varphi}(\lambda+r)\bigg\|
_{\mathcal{B}_{2}}^{2}|\lambda+r|^{n}d\lambda\nonumber\\
&+&\int\limits_{0}^{1}\sum\limits_{r\in\mathbb{Z}}
\sum\limits_{(k_{1}^{\prime},l_{1}^{\prime})\neq (k_{2}^{\prime},l_{2}^{\prime})}\bigg\langle\bigg(\sum\limits_{m^{\prime}}c_{k_{1}^{\prime},l_{1}^{\prime},m^{\prime}}e^{2\pi i m^{\prime}\lambda}\bigg)\widehat{L_{(2k_{1}^{\prime},l_{1}^{\prime},0)}\varphi}(\lambda+r),\nonumber\\
&&\bigg(\sum\limits_{m^{\prime}}c_{k_{2}^{\prime},l_{2}^{\prime},m^{\prime}}e^{2\pi i m^{\prime}\lambda}\bigg)\widehat{L_{(2k_{2}^{\prime},l_{2}^{\prime},0)}\varphi}(\lambda+r)\bigg\rangle
_{\mathcal{B}_{2}}|\lambda+r|^{n}d\lambda,\label{2peq13}
\end{eqnarray}
using \eqref{3peq20}. Now, using \eqref{3peq18} and the fact that $\|\widehat{L_{(2k,l,0)}\varphi}(\lambda)\|_{\mathcal{B}_{2}}=\|\widehat{\varphi}(\lambda)\|_{\mathcal{B}_{2}},~\forall~ k,l\in\mathbb{Z}^{n}$, the first term on the right hand side of \eqref{2peq13} becomes
\begin{eqnarray}
&&\int\limits_{0}^{1}\sum\limits_{r\in\mathbb{Z}}\sum\limits_{k^{\prime},l^{\prime}}
\bigg|\sum\limits_{m^{\prime}}c_{k^{\prime},l^{\prime},m^{\prime}} e^{2\pi im^{\prime}\lambda}\bigg|^{2}\|\widehat{L_{(2k^{\prime},l^{\prime},0)}\varphi}(\lambda+r)\|
_{\mathcal{B}_{2}}^{2}|\lambda+r|^{n}d\lambda\nonumber\\
&=&\int\limits_{0}^{1}\sum\limits_{k^{\prime},l^{\prime}}
\bigg|\sum\limits_{m^{\prime}}c_{k^{\prime},l^{\prime},m^{\prime}} e^{2\pi im^{\prime}\lambda}\bigg|^{2}\sum\limits_{r\in\mathbb{Z}}\|\widehat{\varphi}
(\lambda+r)\|_{\mathcal{B}_{2}}^{2}|\lambda+r|^{n}d\lambda\nonumber\\
&=&\int\limits_{0}^{1}\sum\limits_{k^{\prime},l^{\prime}}
\bigg|\sum\limits_{m^{\prime}}c_{k^{\prime},l^{\prime},m^{\prime}} e^{2\pi im^{\prime}\lambda}\bigg|^{2}G_{0,0}^{\varphi}(\lambda)d\lambda.\label{2peq15}
\end{eqnarray}
The second term on the right hand side of \eqref{2peq13} is
\begin{eqnarray*}
&&\int\limits_{0}^{1}
\sum\limits_{(k_{1}^{\prime},l_{1}^{\prime})\neq (k_{2}^{\prime},l_{2}^{\prime})}\sum\limits_{m_{1}^{\prime},m_{2}^{\prime}}
c_{k_{1}^{\prime},l_{1}^{\prime},m_{1}^{\prime}}\overline{c_{k_{2}^{\prime},l_{2}^{\prime},m_{2}^{\prime}}}e^{2\pi i(m_{1}^{\prime}-m_{2}^{\prime})\lambda}\\
&&\times\sum\limits_{r\in\mathbb{Z}}\langle\widehat{L_{(2k_{1}^{\prime},l_{1}^{\prime},0)}\varphi}(\lambda+r),
\widehat{L_{(2k_{2}^{\prime},l_{2}^{\prime},0)}\varphi}(\lambda+r)\rangle
_{\mathcal{B}_{2}}|\lambda+r|^{n}d\lambda=0,
\end{eqnarray*}
as $\varphi$ satisfies condition C. Thus \eqref{3peq14} follows from \eqref{2peq13} and \eqref{2peq15}. Then, using \eqref{3pr1}, we have
\begin{eqnarray*}
\|f\|^{2}=\sum\limits_{k^{\prime},l^{\prime}}\int\limits_{0}^{1} |\rho_{k^{\prime},l^{\prime}}(\lambda)|^{2} G_{0,0}^{\varphi}(\lambda) d\lambda
&=&\int\limits_{0}^{1} \|\rho(\lambda)\|_{\ell^{2}(\mathbb{Z}^{2n})}^{2}G_{0,0}^{\varphi}(\lambda) d\lambda\\
&=&\|\rho\|_{L^{2}((0,1],\ell^{2}(\mathbb{Z}^{2n}),G_{0,0}^{\varphi})}^{2}.
\end{eqnarray*}
Hence $J(\varphi)$ is an isometry. Using density argument this isometry can be extended to the whole of $V(\varphi)$. Moreover from \eqref{3pr5}, we have $f\in V(\varphi)$ if and only if
\begin{eqnarray}\label{3peq13}
\widehat{f}(\lambda)=\sum\limits_{k^{\prime},l^{\prime}\in\mathbb{Z}^{n}}
\fontdimen16\textfont2=5pt \rho_{k^{\prime},l^{\prime}}(\lambda)\widehat{L_{(2k^{\prime},l^{\prime},
0)}\varphi}(\lambda),
\end{eqnarray}
where $\rho(\lambda)=\{\fontdimen16\textfont2=5pt \rho_{k^{\prime},l^{\prime}}(\lambda)\}_
{(k^{\prime},l^{\prime})\in\mathbb{Z}^{2n}}$ and $\rho\in L^{2}((0,1],\ell^{2}(\mathbb{Z}^{2n}),G_{0,0}^{\varphi})$.
\end{pf}

The following definitions and results are in accordance with \cite{hre}.
\begin{defn}\label{3pdef6}
An element $\varphi\in\mathcal{H}$ is said to be a canonical dual to the system $\{f_{k}:k\in\mathbb{Z}\}$ if $\langle f_{k},\varphi\rangle=\delta_{k,0},~\forall~k\in\mathbb{Z}$. In case, if $f_{k}$ is generated from a single function $f$ by some transformation, then $\varphi$ is called a canonical dual to $f$.
\end{defn}
\begin{defn}
A sequence $\{f_{k}:{k \in\mathbb {Z}}\}$ in $\mathcal{H}$ is said to be Besselian if $\sum\limits_{k\in\mathbb{Z}}c_{k}f_{k}$ is convergent implies $\{c_{k}\}\in\ell^{2}(\mathbb{Z})$.
\end{defn}
\begin{defn}
A sequence $\{f_{k}:{k \in\mathbb {Z}}\}$ in $\mathcal{H}$ is said to be Hilbertian if $\{c_{k}\}\in\ell^{2}(\mathbb{Z})$ implies $\sum\limits_{k\in\mathbb{Z}}c_{k}f_{k}$ converges.
\end{defn}
\begin{lem}\label{3plem1}
Suppose a measurable non-negative function $s$ on $\mathbb{T}^{n}$ satisfies $sm\in L^{1}(\mathbb{T}^{n})$, whenever $m\in L^{1}(\mathbb{T}^{n})$, then $s\in L^{\infty}(\mathbb{T}^{n})$.
\end{lem}
The following theorem is in accordance with Remark 2.3 of \cite{ss}.
\begin{thm}\label{3pth4}
For every measurable subset $A\subset\mathbb{T}^{n}$ with $\lambda{(A)}> 0$, $\lambda{(A^{c})}> 0$, there exists non-zero $f\in L^{2}(\mathbb{T}^{n})$ such that
\begin{enumerate}
\item supp~$f\subset A$;
\item there exists $M> 0$ such that $\|S_{n}(f)\|_{L^{\infty}(\mathbb{T}^{n})}\leq M$ for all $n\in\mathbb{N}$, where $S_{n}(f)$ denotes the nth partial sum of the Fourier series of $f$.
\end{enumerate}
\end{thm}
\section{Canonical dual in twisted shift-invariant spaces}
The following theorem gives a necessary and sufficient condition for the existence of the canonical dual to $\varphi$ in $L^{2}(\mathbb{R}^{2n})$ under condition C. We recall Definition \ref{3pdef6}. Accordingly, we say that a function $\tilde{\varphi}\in L^{2}(\mathbb{R}^{2n})$ is said to be a canonical dual to a function $\varphi$ in $L^{2}(\mathbb{R}^{2n})$ if $\langle T_{(k,l)}^{t}\tilde{\varphi},\varphi\rangle=\delta_{(k,l),(0,0)}$ holds for every $k,l\in\mathbb{Z}^{n}$.
\begin{thm}\label{3pth1}
Let $\varphi\in L^{2}(\mathbb{R}^{2n})$ be such that $\varphi$ satisfies condition C. Then there exists a canonical dual $\tilde{\varphi}$ to $\varphi$ that belongs to $V^{t}(\varphi)$ if and only if $\frac{1}{w_{\varphi}}\in L^{1}(\mathbb{T}^{n})$. Moreover, in this case, $K_{\tilde{\varphi}}(\xi,\eta)=\frac{1}{w_{\varphi}(\xi)} K_{\varphi}(\xi, \eta)$.
\end{thm}
\section{Non-redundancy property of Twisted translates in $L^{2}(\mathbb{R}^{2n})$}
In the following theorem, we shall show that under condition C, the condition $w_{\varphi}(\xi)> 0$ a.e. $\xi\in\mathbb{T}^{n}$ is sufficient for the collection consisting of twisted translates of $\varphi\in L^{2}(\mathbb{R}^{2n})$ to be non-redundant.
\begin{thm}\label{3pth2}
Let $\varphi\in L^{2}(\mathbb{R}^{2n})$ be such that $\varphi$ satisfies condition C. Suppose $w_{\varphi}(\xi)> 0$ a.e. $\xi\in\mathbb{T}^{n}$. Then $\{T_{(k,l)}^{t}\varphi:(k,l)\in \mathbb{Z}^{2n}\}$ is $\ell^{2}$- linearly independent.
\end{thm}
In the following theorem, we shall prove that the converse of the above theorem becomes true without using condition C but under an additional assumption that $\{T_{(k,l)}^{t}\varphi:(k,l)\in \mathbb{Z}^{2n}\}$ is a Bessel sequence.
\begin{thm}\label{3pth3}
Suppose $\{T_{(k,l)}^{t}\varphi:(k,l)\in \mathbb{Z}^{2n}\}$ is a Bessel sequence with bound B that is $\ell^{2}$-linearly independent. Then $w_{\varphi}(\xi)> 0$ a.e. $\xi\in\mathbb{T}^{n}$.
\end{thm}
Now we shall prove the converse of Theorem \ref{3pth2} without assuming that $\{T_{(k,l)}^{t}\varphi:(k,l)\in \mathbb{Z}^{2n}\}$ is a Bessel sequence in Theorem \ref{3pth3}.
\begin{thm}\label{3pth5}
Let $0\neq\varphi\in L^{2}(\mathbb{R}^{2n})$. If $\{T_{(k,l)}^{t}\varphi:(k,l)\in \mathbb{Z}^{2n}\}$ is $\ell^{2}$-linearly independent, then $w_{\varphi}(\xi)> 0$ a.e. $\xi\in\mathbb{T}^{n}$.
\end{thm}
\section{Twisted translates as Besselian and Hilbertian sequences in $L^{2}(\mathbb{R}^{2n})$}
In this section, we shall obtain the characterization for the collection $\{T_{(k,l)}^{t}\varphi:(k,l)\in \mathbb{Z}^{2n}\}$ to be Hilbertian. Towards this, we have the following theorem.
\begin{thm}\label{3pth6}
Let $\varphi\in L^{2}(\mathbb{R}^{2n})$ be such that $\varphi$ satisfies condition C. Assume that $\frac{1}{w_{\varphi}}\in L^{1}(\mathbb{T}^{n})$. If $\{T_{(k,l)}^{t}\varphi:(k,l)\in \mathbb{Z}^{2n}\}$ is Hilbertian, then $||w_{\varphi}||_{\infty}<\infty$.
\end{thm}
Now we shall prove the converse of Theorem \ref{3pth6}.
\begin{thm}\label{3pth7}
Let $\varphi\in L^{2}(\mathbb{R}^{2n})$ be such that $\varphi$ satisfies condition C and $\frac{1}{w_{\varphi}}\in L^{1}(\mathbb{T}^{n})$. Assume that $\|w_{\varphi}\|_{\infty}<\infty$. Then $\{T_{(k,l)}^{t}\varphi:(k,l)\in\mathbb{Z}^{2n}\}$ is Hilbertian.
\end{thm}
Combining Theorem \ref{3pth6}, Theorem \ref{3pth7} and by taking $B=\|w_{\varphi}\|_{\infty}$ in Theorem \ref{3pth10}, we get the following theorem.
\begin{thm}\label{3pth8}
Let $\varphi\in L^{2}(\mathbb{R}^{2n})$ be such that $\varphi$ satisfies condition C and $\frac{1}{w_{\varphi}}\in L^{1}(\mathbb{T}^{n})$. Then the following are equivalent.
\begin{enumerate}
\item [(a)] $\{T_{(k,l)}^{t}\varphi:(k,l)\in\mathbb{Z}^{2n}\}$ is Hilbertian.
\item [(b)] $\|w_{\varphi}\|_{\infty}<\infty$.
\item [(c)] $\{T_{(k,l)}^{t}\varphi:(k,l)\in\mathbb{Z}^{2n}\}$ is a Bessel sequence.
\end{enumerate}
\end{thm}
\begin{thm}\label{3pth9}
Let $\varphi\in L^{2}(\mathbb{R}^{2n})$ be such that $\varphi$ satisfies condition C and $\frac{1}{w_{\varphi}}\in L^{1}(\mathbb{T}^{n})$. If any one of the equivalent conditions of Theorem \ref{3pth8} is true, then there exists $\tilde{\varphi}\in V^{t}(\varphi)$, canonical dual to $\varphi$ such that $\{T_{(k,l)}^{t}\tilde{\varphi}:(k,l)\in\mathbb{Z}^{2n}\}$ is Besselian.

Moreover, if $\{T_{(k,l)}^{t}\varphi:(k,l)\in\mathbb{Z}^{2n}\}$ is Besselian and if for each $\psi\in V^{t}(\varphi)$ there exists a sequence $\{c_{k,l}\}$ such that $\psi=\sum\limits_{k,l\in\mathbb{Z}^{n}}c_{k,l}T_{(k,l)}^{t}\varphi$, then $\|\frac{1}{w_{\varphi}}\|_{\infty}<\infty$.
\end{thm}
\section{Results on the Heisenberg group}
The following theorem provides a necessary and sufficient condition for the existence of the canonical dual to a function $\varphi\in L^{2}(\mathbb{H}^{n})$ under condition C. We recall Definition \ref{3pdef6}. Correspondingly, we say that a function $\tilde{\varphi}\in L^{2}(\mathbb{H}^{n})$ is said to be a canonical dual to a function $\varphi$ in $L^{2}(\mathbb{H}^{n})$ if $\langle L_{(2k,l,m)}\tilde{\varphi},\varphi\rangle=\delta_{(k,l,m),(0,0,0)}$ holds for every $k,l\in\mathbb{Z}^{n}, m\in\mathbb{Z}$.
\begin{thm}\label{3pth12}
Let $\varphi\in L^{2}(\mathbb{H}^{n})$ be such that $\varphi$ satisfies condition C. Then there exists a canonical dual $\tilde{\varphi}$ to $\varphi$ that belongs to $V(\varphi)$ if and only if $\frac{1}{G_{0,0}^{\varphi}}\in L^{1}(0,1]$. Moreover, in this case, $\widehat{\tilde{\varphi}}(\lambda)=\frac{1}{G_{0,0}^{\varphi}(\lambda)}
\widehat{\varphi}(\lambda)$.
\end{thm}
\begin{thm}
Let $\varphi\in L^{2}(\mathbb{H}^{n})$ be such that $\varphi$ satisfies condition C. If $G_{0,0}^{\varphi}(\lambda)>0$ a.e. $\lambda\in (0,1]$, then $\{L_{(2k,l,m)}\varphi:(k,l,m)\in\mathbb{Z}^{2n+1}\}$ is $\ell^{2}$-linearly independent.
\end{thm}
In the following, we state the converse of the above theorem without proof. The proof will be similar to the proof of Theorem \ref{3pth5}.
\begin{thm}
Let $0\neq\varphi\in L^{2}(\mathbb{H}^{n})$. Assume that $\{L_{(2k,l,m)}\varphi:(k,l,m)\in\mathbb{Z}^{2n+1}\}$ is $\ell^{2}$-linearly independent. Then $G_{0,0}^{\varphi}(\lambda)>0$ a.e. $\lambda\in (0,1]$.
\end{thm}
As in the case of twisted translation, by taking $B=\|G_{0,0}^{\varphi}\|_{\infty}$ in Theorem \ref{3pth11}, we can show the following theorem.
\begin{thm}
Let $\varphi\in L^{2}(\mathbb{H}^{n})$ be such that $\varphi$ satisfies condition C and $\frac{1}{G_{0,0}^{\varphi}}\in L^{1}(0,1]$. Then the following are equivalent.
\begin{enumerate}
\item [(a)] $\{L_{(2k,l,m)}\varphi:(k,l,m)\in\mathbb{Z}^{2n+1}\}$ is Hilbertian.
\item [(b)] $\|G_{0,0}^{\varphi}\|_{\infty}<\infty$.
\item [(c)] $\{L_{(2k,l,m)}\varphi:(k,l,m)\in\mathbb{Z}^{2n+1}\}$ is a Bessel sequence.
\end{enumerate}
\end{thm}
\begin{thm}
Let $\varphi\in L^{2}(\mathbb{H}^{n})$ be such that $\varphi$ satisfies condition C and $\frac{1}{G_{0,0}^{\varphi}}\in L^{1}(0,1]$. If any of the equivalent conditions of the above theorem holds, then there exists $\tilde{\varphi}\in V(\varphi)$ such that $\{L_{(2k,l,m)}\tilde{\varphi}:(k,l,m)\in\mathbb{Z}^{2n+1}\}$ is Besselian.

Moreover, if $\{L_{(2k,l,m)}\varphi:(k,l,m)\in\mathbb{Z}^{2n+1}\}$ is Besselian and if for each $\psi\in V(\varphi)$, there exists a sequence $\{c_{k,l,m}\}$ such that $\psi=\sum\limits_{(k,l,m)\in\mathbb{Z}^{2n+1}}c_{k,l,m}L_{(2k,l,m)}\varphi$, then $\|\frac{1}{G_{0,0}^{\varphi}}\|_{\infty}<\infty$.
\end{thm}


\begin{thebibliography}{99}
\bibitem{bhm} Barbieri, D., Hernandez, E. and Mayeli, A.: Bracket map for the Heisenberg group and the characterization of cyclic subspaces. Appl. Comput. Harmon. Anal. 37 (2014) 218-234.
\bibitem{bhp} Barbieri, D., Hernandez, E. and Parcet, J.: Riesz and frame systems genarated by unitary actions of discrete groups. Appl. Comput. Harmon. Anal. In Press (2015).
\bibitem{bo} Bownik, M.: The structure of shift-invariant Subspaces of $L^{2}(\mathbb{R}^{n})$. J. Funct. Anal. 176 (2000), 282-309.
\bibitem{cab} Cabrelli, C. and Paternostro, V.: Shift-invariant spaces on LCA groups. J. Funct. Anal. 258 (2010) 2034-2059.
\bibitem{chr2} Christensen, O.: Frames and Bases, An introductory course. Birkh$\ddot{a}$user, Boston, 2008.
\bibitem{cmo} Currey, B., Mayeli, A. and Oussa, V.: Characterization of shift-invariant spaces on a class of Nilpotent Lie groups with application. J. Fourier Anal. Appl. 20, 384-400, 2014.
\bibitem{heil} Heil, Christopher.: A basis theory primer. Birkh$\ddot{a}$user, Expanded edition, 2011.
\bibitem{hre} Hernandez, E., Sikic, H., Weiss, G. and Wilson, E.: On the properties of the integer translates of a square integrable function in $L^{2}(\mathbb {R})$, Harmonic analysis and partial differential equations, Contemp. Math. Vol. 505, Amer. Math. Soc. Providence, RI, 2010, PP (233-249).
\bibitem{kgr}  Kamyabi Gol, R.A. and Raisi Tousi, R.: The structure of shift invariant spaces on a locally compact abelian group.
    J. Math. Anal. Appl. 340, 219-225, 2008.
\bibitem {ras} Radha, R. and Adhikari, S.: Frames and Riesz bases of twisted shift-invariant spaces in $L^{2}(\mathbb{R}^{2n})$. J. Math. Anal. Appl. 434(2016) 1442-1461.
\bibitem {ra} Radha, R. and Adhikari, S.: Shift-invariant spaces with countably many mutually orthogonal generators on the Heisenberg group. arxiv:1711.06902 [math.FA].
\bibitem{rad} Radha,R. and Shravan kumar, N.: Shift-invariant subspaces on compact groups, Bull. Sci. Math. Vol 137, No 4, 485-497, 2013.
\bibitem{lf} Luef, F.: Projections in noncommutative tori and Gabor frames. Proc. Amer. Math. Soc. 139 (2011), 571-582.
\bibitem{ss} Saliani, S.: $\ell^{2}$-linear independence for the system of integer translates of a square integrable function. Proc. Amer. Math. Soc. Volume 141, No 3, 937-941, 2013.
\bibitem{th} Thangavelu, S.: Harmonic analysis on the Heisenberg group. Birkh$\ddot{a}$user, Boston, 1997.
\end{thebibliography}
\end{document}